\documentclass[journal,12pt,onecolumn,draftclsnofoot,]{ieeeconf} 

\IEEEoverridecommandlockouts                              

\usepackage{graphicx}      
\graphicspath{{./figures/}}
\usepackage{xcolor}
\usepackage[utf8]{inputenc}

\usepackage{subfigure}

\usepackage{hyperref}

\usepackage{amsmath,amssymb,amsfonts}
\usepackage{algorithmic}
\usepackage{textcomp}
\usepackage{mathtools}
\usepackage{placeins}
\usepackage{siunitx}
\usepackage{bm}
\newcommand{\bs}[1]{\boldsymbol{{#1}}}

\makeatletter
\newcommand\deflabel[1]{\def\@currentlabel{#1}}

\makeatother

\title{\LARGE \bf
Spatio-Temporal Decomposition of Sum-of-Squares Programs for the Region of Attraction and Reachability
}

\author{Vít Cibulka$^{1,2}$, Milan Korda$^{1,2}$ and Tomáš Haniš$^{1}$
\thanks{*This research was supported by the Czech Science Foundation (GACR)
 under contracts No. GA19-18424S, GA20-11626Y, 
 and by the Grant Agency of the Czech Technical University in Prague,
grant No. SGS19/174/OHK3/3T/13. This work has also been supported by European Union’s Horizon 2020 research and innovation programme under the Marie Skłodowska-Curie Actions, grant agreement 813211 (POEMA) and by the AI Interdisciplinary Institute ANITI funding, through the
French “Investing for the Future PIA3” program under the Grant agreement n$^\circ$ ANR-19-PI3A-0004. 
}
\thanks{$^{1}$Department of Control Engineering, Faculty of Electrical Engineering,
Czech Technical University in Prague, The Czech Republic
        {\tt\small vit.cibulka@fel.cvut.cz, tomas.hanis@fel.cvut.cz}}%
\thanks{$^{2}$ CNRS, Laboratory for Analysis and Architecture of Systems, Toulouse, France
        {\tt\small korda@laas.fr}}%
		\thanks{This work has been submitted to the IEEE for possible publication. Copyright may be transferred without notice, after which this version may no longer be accessible.}
}

\begin{document}

\maketitle
\thispagestyle{empty}
\pagestyle{empty}

\begin{abstract}
This paper presents a method for calculating \mbox{Region of Attraction} of a target set (not necessarily an equilibrium)
 for controlled polynomial dynamical systems, using a hierarchy of semidefinite programming problems (SDPs).
 Our approach builds on previous work and addresses its main issue, the fast-growing memory demands for solving large-scale SDPs.

The main idea in this work is in dissecting the original 
resource-demanding problem into multiple smaller, interconnected, and easier to solve problems.
 This is achieved by spatio-temporal splitting akin to methods based on partial differential equations.
We show that the splitting procedure retains the convergence and outer-approximation guarantees of the previous work, while  achieving higher precision in less time and with smaller memory footprint.

\end{abstract}


\section{Introduction}

This paper deals with stability and reachability analysis of nonlinear dynamical systems.  
A prominent method for quantifying the stability of a nonlinear system is based on calculating 
its Region Of Attraction (ROA) with respect to a given target set, which is the topic of this work.
 Using a time-reversal, one obtains the reachable set. Both of these objects are crucial in assessing stability and safety of control system; as a concrete example, let us mention the Flight Envelopes in 
the aerospace industry that directly coincide with the ROA~\cite{roa_flight_envelope}.
 A more recent application would be the study of closed-loop parameters variation, of 
 both the controller and the controlled plant, and its effect on the ROA.
For example, given a vehicle on the road, these  parameters could be the distribution of mass of the vehicle or various road conditions~\cite{roa_car}.

The dominant methods for computing ROA are based on Lyapunov functions whose level sets provide inner approximations the ROA \cite{Khalil2002}. For polynomial systems, the Lyapunov function can be found 
by solving semidefinite programming problems (SDPs) \cite{Chesi2004} and \cite{Tibkena}. These methods only work for autonomous systems and ROA's with respect to given equilibria.

This work expands the approach presented in \cite{milanroa}, which is based on optimizing
over trajectories of a polynomial dynamical systems modelled as occupation measures and is not limited to an autonomous 
system nor a stable equilibrium.  The contribution proposed in this paper is dissection of the original problem into multiple interconnected problems of lower complexity. The original formulation is discretized in both the time and state variables, aiming to strike a trade-off between the original optimization-based approach~\cite{milanroa} where a single polynomial is defined on the entire state-space and various methods based on extensive discretization such as the PDE-based approach \cite{Mitchell2003}, set oriented methods \cite{Dellnitz2001} or transfer operator approaches \cite{Wanga}. By doing so, our approach, based fully on convex optimization, preserves the outer approximation and convergence guarantees of \cite{milanroa}, both of which are explicitly proven in this work, while significantly improving scalability of the method.

\textbf{Structure of this paper}
Section \ref{sec:problem_statement} presents the problem statement and
Section \ref{sec:splitting} introduces the splitting procedure with the proofs of outer approximation and convergence.
Section \ref{sec:sos} states the practical sum-of-squares variant of the 
problem which is then demonstrated in Section \ref{sec:numerical_examples} on 
numerical examples.
The paper ends with a conclusion and future work discussion in sections 
\ref{sec:conclusion} and \ref{sec:future_work} respectively.

\textbf{Notation}
The symbol \(\mathbb{Z}_{j}\) denotes the set of consecutive 
integers \(\{1,2,\dots,j\}\). The Lebesgue measure (i.e., the volume) of a set is denoted by $\lambda(A)$. The indicator function \(I_{A}(x)\) of a set $A$ is the function that takes the value \(1\) for \(x \in A\), and \(0\) otherwise. The interior of a set $A$ is denoted by $A^\circ$.
The symbols \(C(A)\) and \(C^1(A)\) denote respectively the sets of continuous and continuously differentiable functions on \(A\). The set of Borel measurable functions defined on a set $A$ taking values in a set $B$ is denoted by $L(A;B)$ 

\section{Problem statement}
\label{sec:problem_statement}
Let us consider the nonlinear system with control
\begin{equation}
	\label{eq:nlsys}
	\dot{x}(t) = f(t,x(t),u(t)), t \in [0,T],
\end{equation}
where \(x(t) \in \mathbb{R}^n\) is the state vector, \(u(t) \in \mathbb{R}^m\)
is the control input vector, \(t\) is time, \(T > 0\) is the final time 
and \(f\) is the vector field, which is assumed to be polynomial in variables \(x\) and \(u\).

The state and control input are constrained by basic semialgebraic sets
\begin{equation}
	\label{eq:semialg_sets}
\begin{alignedat}{2}
 	&u(t) \in U   :=&&\{u \in \mathbb{R}^m : g_j^U(u)       \geq 0, j \in \mathbb{Z}_{n_U}\}, t\in [0,T] ,\\
 	&x(t) \in X   :=&&\{x \in \mathbb{R}^n : g_j^X(x)       \geq 0, j \in \mathbb{Z}_{n_X}\}, t\in [0,T] ,\\
 	&x(T) \in X_T :=&&\{x \in \mathbb{R}^n : g_j^{ X_T }(x) \geq 0, j \in \mathbb{Z}_{N_{X_T}}\},
\end{alignedat}
\end{equation}
where \(g_j^U(u)\),  \(g_j^X(x) \), and \(g_j^{X_T}(x)\) are polynomials. The region of attraction (ROA) is then defined as 
\begin{equation}
\begin{alignedat}{2}
	X_0 = \{& x_0 \in &&X : \exists\, u(\cdot) \in L([0,T];U) \\
	&\text{s.t.}&& \dot{x} = f(t,x(t),u(t)) \text{ a.e. on}\;  [0,T],\\
	& &&x(0) = x_0,\;x(t) \in X\;\forall\, t\in [0,T],\; x(T) \in X_T  \},
\end{alignedat}
\end{equation}
where ``a.e.'' stands for ``almost everywhere'' with respect to the Lebesgue measure.




\section{Time and state space splitting}
\label{sec:splitting}

It was shown in \cite{milanroa}, that the ROA can be characterized by an infinite-dimensional linear programming (LP) problem in the space of Borel measures or by its dual in the space of continuous functions. These LPs can then be approximated by SDPs, with guaranteed convergence. For space reason, the entire exposition in this work will be in the ``dual'' setting of continuous functions approximated by polynomials.

Let us restate the original problem from \cite{milanroa}:
\begin{equation}
	\label{eq:dualOriginal}
	\begin{alignedat}{2}
		d^\star &= \inf \int_{X} w(x)d\lambda(x)             &&,\\
		\text{s.t.}& ( \mathcal{L}v )(t,x,u) \leq 0         &&,\forall (t,x,u) \in [0,T] \times X \times U\\
		&w(x)                \geq v(0,x) + 1&&,\forall x \in X\\
		&v(T,x)              \geq 0         &&,\forall x \in X_T\\
		&w(x)                \geq 0         &&, \forall x \in X,\\
	\end{alignedat}
\end{equation}

with variables \(w(x) \in C(X)\) and \(v(t,x) \in C([0,T]\times X)\).

Any minimizing sequence \((w_k,v_k)\) for (\ref{eq:dualOriginal}) satisfies \( w_k \ge I_{X_0}(x)\) and $w_k \to I_{X_0}$ in $L_1$  as well as \mbox{\(\{x \in X : v_k(0,x) \geq 0\} \supset X_0\)} with convergence in terms of the volume discrepancy tending to zero (see \cite{milanroa} for proof).

Let us now split the state space \(X\) into \(I\) closed subsets~\(X_i\) 
\begin{equation}
	X = \bigcup\limits_{i=1}^{I} X_{i}
\end{equation}
and the time interval \([0,T]\) into \(K-1\) intervals \([T_k, T_{k+1}]\) 
\begin{equation}
	[0,T] = \bigcup\limits_{k=1}^{K-1} [T_k,T_{k+1}],
\end{equation}
where  \(K\) is the number of time splits (meaning \(K-1\) intervals). It is assumed that $X_i^\circ \cap X_j^\circ = \emptyset$ for $i \ne j$.

The function \(w(x)\) will be split into \(I\) functions \(w_i(x)\) 
\begin{equation}
\label{eq:split_w}
	w(x) = 
	\begin{cases}
		w_1(x) &\text{for } x \in X_1\\
		\dots \\
		w_i(x) &\text{for } x \in X_i\\
		\dots\\
		w_I(x) &\text{for } x \in X_I
	\end{cases}
\end{equation}
and \(v(t,x)\) will be split into \(I\cdotp(K-1)\) functions \(v_{i,k}(t,x)\)
\begin{equation}
	\label{eq:split_v}
	v(t,x) = 
	\begin{cases}
		v_{1,1}(t_1,x_1) &\text{for } t \in [T_1,T_{2}], x \in X_1\\
		\dots\\
		v_{i,k}(t_k,x_i) &\text{for } t \in [T_k,T_{k+1}], x \in X_i\\
		\dots\\
		v_{I,K-1}(t_{K-1},x_I) &\text{for } t \in [T_{K-1},T_{K}], x \in X_I.
	\end{cases}
\end{equation}

Assuming that neighbouring subsets \(X_i\) share boundaries,
let us define the set of indices of these neighbours as
\begin{equation}
\begin{alignedat}{2}
N_X := \{(a,b): X_a \cap X_b \not= \emptyset\}.
\end{alignedat}
\end{equation}

We are now ready to write the split version of \eqref{eq:dualOriginal} 
\begin{equation}
	\label{eq:dualOriginalSPlit}
	\begin{split}
	\begin{alignedat}{3}
		& d_\mathrm{s}^\star = &&\inf \sum_i \int_{X_i} w_i(x)d\lambda(x)     \\
& \text{s.t.\,}&& \text{for all } i \in \mathbb{Z}_{I}, k \in \mathbb{Z}_{K-1} \text{ and }\\
& && \hspace{3cm} (a,b) \in N_X, &&x_{a,b}\in X_a\cap X_b\\
		&&& ( \mathcal{L}v_{i,k} )(t,x,u) \leq 0  \quad\,  \forall(t,x,u)  \in [&&T_k,T_{k+1}] \times X_i \times U    \\
	&	&&w_i(x)                \geq v_{i,k}(0,x) + 1  \quad &&\forall x \in X_i   \\
&		&&v_{i,K}(T,x)              \geq 0   \quad &&\forall x \in X_T \\
&    &&w_i(x)                \geq 0         \quad &&\forall x \in X_i \\
&	&& v_{i,k}(T_{k+1},x) \geq v_{i,k+1}(T_{k+1},x) \quad &&\forall x \in X_i\\
&	&&( v_{a,k}(t,x_{a,b}) - v_{b,k}(t,x_{a,b}) )\cdot  h_{a,b}^\top&&f(t,x_{a,b},u)  \geq 0
	\end{alignedat}\\
\end{split}
\end{equation}
 where \(h_{a,b}^\top\) is a normal vector of a shared boundary between 
two neighbouring sets \(X_a\) and \(X_b\), and \(x_{a,b}\) is a set of points on said
boundary such that \( x_{a,b} \in X_a \cap X_b\). For simplicity we assume that the normal vector $h_{a,b}$ is independent of $x$; the case of polynomial or rational dependence of $h$ on $x$ can also be handled~\cite[Section 4.3]{korda2018moments}. The optimization variables in~\eqref{eq:dualOriginalSPlit} are the continuously differentiable functions of $v_{i,k}$, each defined some neighborhood of $[T_k,T_{k+1}]\times X_i$ and the continuous functions $w_i$, each defined on $X_i$. 

Let us now show that the modified problem \eqref{eq:dualOriginalSPlit} provides a guaranteed outer approximation of \(X_0\) which can be defined as
\begin{equation}
	\label{eq:def_X0v}
	\bar{X}_{v,0} := \{x : v(0,x) \geq 0\}.
\end{equation}
\newtheorem{theorem}{Theorem}
\numberwithin{paragraph}{theorem}
\begin{theorem}
	For any pair \((v,w)\) feasible in \eqref{eq:dualOriginalSPlit}, 
	it holds that \(v(0,\cdot) \geq 0\) on \(X_0\) and $\bar{X}_{v,0} \supset X_0$.
\end{theorem}
\begin{proof}
We first need to show that the discontinuous function \(v(t,x(t))\) is 
decreasing along the system trajectories. That is, given two time instants \(t_\alpha \leq t_\beta\),  we want to show that 
	\begin{equation}
		\label{eq:cont_flow_cond}
	v(t_\beta,x(t_\beta)) \leq v(t_\alpha,x(t_\alpha)).
	\end{equation}

If we $v$ were to be differentiable, this follows by simply integrating the first constraint of (\ref{eq:dualOriginal}) along a trajectory.

	Let us now show that \eqref{eq:cont_flow_cond} holds even for 
	the discontinuous \(v(t,x(t))\) as defined in \eqref{eq:split_v}. The first constraint of \eqref{eq:dualOriginalSPlit} ensures that the values of $v$ decrease whenever the trajectory resides in the interior of one of the sets $X_i$. We therefore need to argue only about what happens on the boundary of these sets.

    \paragraph{Time splits}
	\label{sec:subproof_time}
	The result for time splits follows immediately from 	the fifth constraint of~(\ref{eq:dualOriginalSPlit}).
	
%

\paragraph{State-space splits}
\label{sec:subproof_state}
Let us first assume that the state space \(X\) is split into two parts,
\(X_\alpha\) and \(X_\beta\) by a hyper-plane with normal vector 
 \(h\), pointing from \(X_\alpha \) to \(X_\beta\), so that 
\begin{equation}
\begin{alignedat}{2}
h^\top (x_\beta - x_\alpha) \geq 0,
\end{alignedat}
\end{equation}
for \(x_\alpha \in X_\alpha\) and \(x_\beta \in X_\beta\).

The function \(v(t,x(t))\), now split between \(X_\alpha\) and \(X_\beta\),
is defined as
\begin{equation}
	\label{eq:v_def_Xsplits}
	v(t,x(t)) = 
	\begin{cases}
		v_\alpha(t,x(t)) &\text{for } t \in [0,T], x(t) \in X_\alpha\\
		v_\beta(t, x(t)) &\text{for } t \in [0,T], x(t) \in X_\beta.
	\end{cases}
\end{equation}

Let $x_0 \in X_0 \cap X_\alpha$, $t_\alpha \in [0,T]$ and $u(\cdot)$ be given. Let \(x(\cdot|x_0)\) be the trajectory starting at $t_\alpha$ generated by $u(\cdot)$ and suppose that  $x(t\mid x_0) \in X$, $u(t) \in U$ for $t \in [t_\alpha,T]$. Assume further that this trajectory crosses from \(X_\alpha\) to \(X_\beta\) at the crossing time
\begin{equation}
	\tau = \inf_{t} \{t_\alpha \leq t \mid  x(t | x_0) \in X^\circ_\beta\} \le T
\end{equation}
and assume that this trajectory stays in $X_\beta$ for $t \in [\tau,T]$. At the crossing point $x(\tau)$, it holds
\begin{equation*}
	h^\top f(\tau,x(\tau),u(\tau)) \geq 0.
\end{equation*}
The last constraint of \eqref{eq:dualOriginalSPlit} implies that
\begin{equation}
	\label{eq:flow_cond}
 v_\alpha(\tau,x(\tau)) \geq v_\beta(\tau,x(\tau))
\end{equation}
whereas the first constraint implies
\begin{equation}\label{eq:proofAux1}
\frac{\mathrm{d}}{\mathrm{d}t} v_\alpha(t,x(t)) \leq 0, t\in[t_\alpha,\tau)
\end{equation}
and
\begin{equation}\label{eq:proofAux2}
\frac{\mathrm{d}}{\mathrm{d}t} v_\beta(t,x(t)) \leq 0, t\in(\tau,t_\beta].
\end{equation}


Let us now calculate the value of \(v_\beta(t_\beta, x(t_\beta))\):
\begin{multline}
v_\beta(t_\beta, x(t_\beta)) =v_\alpha(t_\alpha, x(t_\alpha)) + \int_{t_\alpha}^\tau \frac{\mathrm{d}}{\mathrm{d}t}v_\alpha(t,x(t)) \mathrm{d}t\\
+ [v_\beta(\tau, x(\tau)) - v_\alpha(\tau, x(\tau))] +  \int_{\tau}^{t_\beta} \frac{\mathrm{d}}{\mathrm{d}t}v_\beta(t,x(t)) \mathrm{d}t.
\end{multline}
By inspecting \eqref{eq:flow_cond}, \eqref{eq:proofAux1} and \eqref{eq:proofAux1}, we can see that
the last three summands are nonpositive and we get the inequality $v_\beta(t_\beta, x(t_\beta)) \leq v_\alpha(t_\alpha, x(t_\alpha))$, which is equivalent to \eqref{eq:cont_flow_cond},  recalling the definition of \(v(t,x(t))\) in \eqref{eq:v_def_Xsplits}. The procedure for the negative trajectory direction is analogous. We note that this analysis encompasses the subtle case of the trajectory sliding on the boundary between the sets $X_\alpha$ and $X_\beta$.

By induction, we can prove the inequality for arbitrary splitting of the state-space and 
time axis by considering a sequence of crossing times associated to a given trajectory. Therefore $v(t_\beta,x(t_\beta)) \leq v(t_\alpha,x(t_\alpha))$ for any $0\le t_\alpha \le t_\beta \le T$. By setting \(t_\alpha = 0, t_\beta = T\), and using the constraints 
of \eqref{eq:dualOriginalSPlit}, we get
\begin{equation}
\begin{alignedat}{2}
	v(T,x(T)) &\leq v(0,x_0)\\
	0 \leq v(T,x(T)) &\leq v(0,x_0)\\
	0 &\leq v(0,x_0)
\end{alignedat}
\end{equation}
for any $x_0 \in X_0$ as desired. This also implies that $x_0 \in \bar X_{v,0}$ and hence $\bar X_{v,0} \supset X_0$.

%

%

\end{proof}


\section{SOS representation}
\label{sec:sos}
We can now obtain the SDP representation of \eqref{eq:dualOriginalSPlit} by applying
Putinar's Positivstellensatz~\cite{putinar}.
For example, given polynomials \(c(x)\) and \(g(x)\) the inequality 
\begin{equation}
	c(x) \geq 0 \quad \text{for } x \in \{x: g(x) \geq 0\}
\end{equation}
is implied by 
\begin{equation}
	c(x)  = q(x) + s(x)g(x),
\end{equation}
where \(q(x)\) and \(s(x)\) are sum-of-squares polynomials. The condition that a polynomial $s$ of degree $2d$ is sum-of-squares is in turn equivalent to $s(x) = m(x)^\top W m(x)$, $W \succeq 0$, where $m(x)$ is a basis of polynomials up degree $d$ and hence this constraint is SDP representable.

The SOS approximation  of \eqref{eq:dualOriginalSPlit} reads
\begin{equation}
	\label{eq:dualSOS}
	\begin{alignedat}{3}
		& &&  \inf \sum_i \text{ }  \mathbf{w}_i^\top l_i\\
&\text{s.t.\quad} &&\text{for all } i \in \mathbb{Z}_{I}, k \in \mathbb{Z}_{K-1} \text{ and }  (a,b) \in N_X \\
		& &&-(\mathcal{L}v_{i,k})(z) =  q_{i,k}(z) + \bs{s^\tau_{i,k}}(z)^\top\bs{g^{\tau}_k}(t)\\
		& &&\hspace{2.2cm} + \bs{s^X_{i,k}}(z)^\top \bs{g^X_{i,k}}(x) + \bs{s^U_{i,k}}(z)^\top \bs{g^U}(u) &&\\
		& &&w_i(x) - v_{i,k}(0,x) - 1 = q_{0_{i,k}}(x) + \bs{s^0_{i,k}}(x)^\top \bs{g_{i}^X}(x)\\
		& &&v_{i,K}(T,x) = q^T_i(x) + \bs{s^{X_T}_{i}}(x)^\top \bs{g^{X_T}_{i}}(x)\\
		& &&w_i(x) = q^w_i(x) + \bs{s^w_{i}}(x)^\top \bs{g^X_{i}}(x) \\
		& &&v_{i,k}(T_{k+1},x) - v_{i,k+1}(T_{k+1},x) = \\
		& &&\hspace{3.35cm}q_{{i,k}}^\tau(x) + \bs{s^t_{i,k}}(x)^\top \bs{g^X_{i}}(x) \\
		& && ( v_{a,k}(t,x_{a,b}) - v_{b,k}(t,x_{a,b}) ) = \\
		& &&\hspace{0.59cm} q^{1}_{k,a,b}(z) + \sum\nolimits_{j=1}^{n_X} s^{1}_{j,k,a,b}(z) h_{a,b}^\top f(t,x_{a,b},u)&& \\
		& && ( v_{b,k}(t,x_{a,b}) - v_{a,k}(t,x_{a,b})  ) = \\
		& &&\hspace{0.59cm} q^{2}_{k,a,b}(z) - \sum\nolimits_{j=1}^{n_X} s^{2}_{j,k,a,b}(z) h_{a,b}^\top f(t,x_{a,b},u)&&
	\end{alignedat}
\end{equation}
where \(z = [t,x,u]^\top\),
\(w_i(x)\) and \(v_{i,k}(t,x)\) are polynomials,
\(\mathbf{w}_i\) is a vector of coefficients of \(w_i(x)\) and \(l_i\) is a 
vector of Lebesgue measure 
moments indexed with respect to the same basis as the coefficients of $w_i$.  The decision variables in the problem are the polynomials $v_{i,k}$ and $w_i$ as well as sum-of-squares multipliers $q$, $s$ and $\bs s$.  The symbols $\boldsymbol{g_{i}^X}$, $\boldsymbol{g_{i}^{X_T}}$, $\boldsymbol{g_{i}^U}$ and $\boldsymbol{g_{k}^\tau}$ denote the column vectors of polynomials describing the sets $X_i$, $X_T\cap X_i$,  $U$ and $[T_k,T_{k+1}]$ in that order.  The degrees of all polynomial decision variables is chosen such that the degrees of all polynomials appearing in (\ref{eq:dualSOS}) do not exceed a given relaxation order $d$. This is a design parameter controlling the accuracy of the approximation.

Given the picewise polynomial functions $(w^d,v^d)$ of degree $d$ constructed
 from a solution to \eqref{eq:dualSOS} as in (\ref{eq:split_w}) and (\ref{eq:split_v}), 
 the outer approximation to the ROA is defined by
\[
X_d = \{x \mid v^d(0,x) \ge 0\}.
\]

Convergence of the SDP approximations holds under the classical 
Archimedianity assumption, e.g., \cite[Assumption 3]{milanroa}.

\begin{theorem}
For each $d \in \mathbb{N}$, we have $X_d \supset X_0$. If in addition the algebraic description of each element of the space-time partition used in~(\ref{eq:dualSOS}) satisfies the Archimedianity condition, then $\lim_{d\to\infty} \lambda(X_d \setminus X_0) = 0$.
\end{theorem}
\begin{proof}
The result follows from \cite[Theorem 6]{milanroa} by taking all functions defined on the partition equal, i.e.,  $v_{i,k} = v$ and $w_i = w$ for some polynomials $v$ and $w$;  this leads to the setting of \cite[Theorem 6]{milanroa}.
\end{proof}

\subsection{Practical implications}
\label{sec:practical_impl}
The ROA with splits is expected to improve accuracy of the original formulation by allowing one to trade off the degree of the polynomials for number of splits. By increasing the degree \(d\), the size of the SDP will increase
with the rate of the binomial coefficient \({n+d}\choose{n}\).
By fixing \(d\) and splitting the state-space into \(\lambda\) cells, 
the SDP size will  grow only linearly with the number of the cells with rate \(\lambda {n+d\choose{n}}\). The numerical results in the following section suggest that even the computation time  grows linearly. A full quantitative understanding of this tradeoff is beyond the scope of the present work; the first step in this direction would be the convergence rate bounds of~\cite{korda2017convergence}. 
\section{Numerical examples}
\label{sec:numerical_examples}
This section presents numerical examples, showcasing the 
performance difference between the proposed method and the 
original approach from \cite{milanroa}.

The first example in Section~\ref{sec:univar} shows the influence of the split positions on the 
resulting ROA.
It is shown that by having the splits exactly at the boundaries of the ROA, we can retrieve
the theoretical indicator function \(I_{X_0}(x)\).

The second example Section~\ref{sec:brockett} benchmarks the algorithm on a Brockett integrator,
which mimics a kinematic model of a nonholonomic system 
(it can be shown that three-dimensional \mbox{nonholonomic} vehicle
with two inputs can be transformed into the 
Brockett integrator \cite{brocket_transform}).


Finally, Section~\ref{sec:performance} presents a comparison of computational demands of the 
proposed method and the original one from \cite{milanroa}.

All the examples were implemented in MATLAB \cite{MATLAB:2020} with 
the use of YALMIP \cite{yalmip}.
The YALMIP's sum-of-squares package \cite{yalmip_sos} was used for rapid
prototyping; for larger examples, the SDPs were assembled using a custom routine. All SDP's were solved by MOSEK \cite{mosek}.

\subsection{Univariate cubic dynamics}
\label{sec:univar}
This example shows that one can find the ROA with a very low degree 
polynomials by correctly positioning the splits.

The system in question is defined as
\begin{equation}
	\dot{x} = x(x - 0.5)(x + 0.5)
\end{equation}
with the state space 
\(X = [-1,1]\), the target set \mbox{\(X_T = [-0.01, 0.01]\)} and 
terminal time \(T = 100\).
The analytic solution of the ROA is \(X_0 = [-0.5, 0.5]\).

\begin{figure}[]
	\centering
	\includegraphics[width=0.5\textwidth]{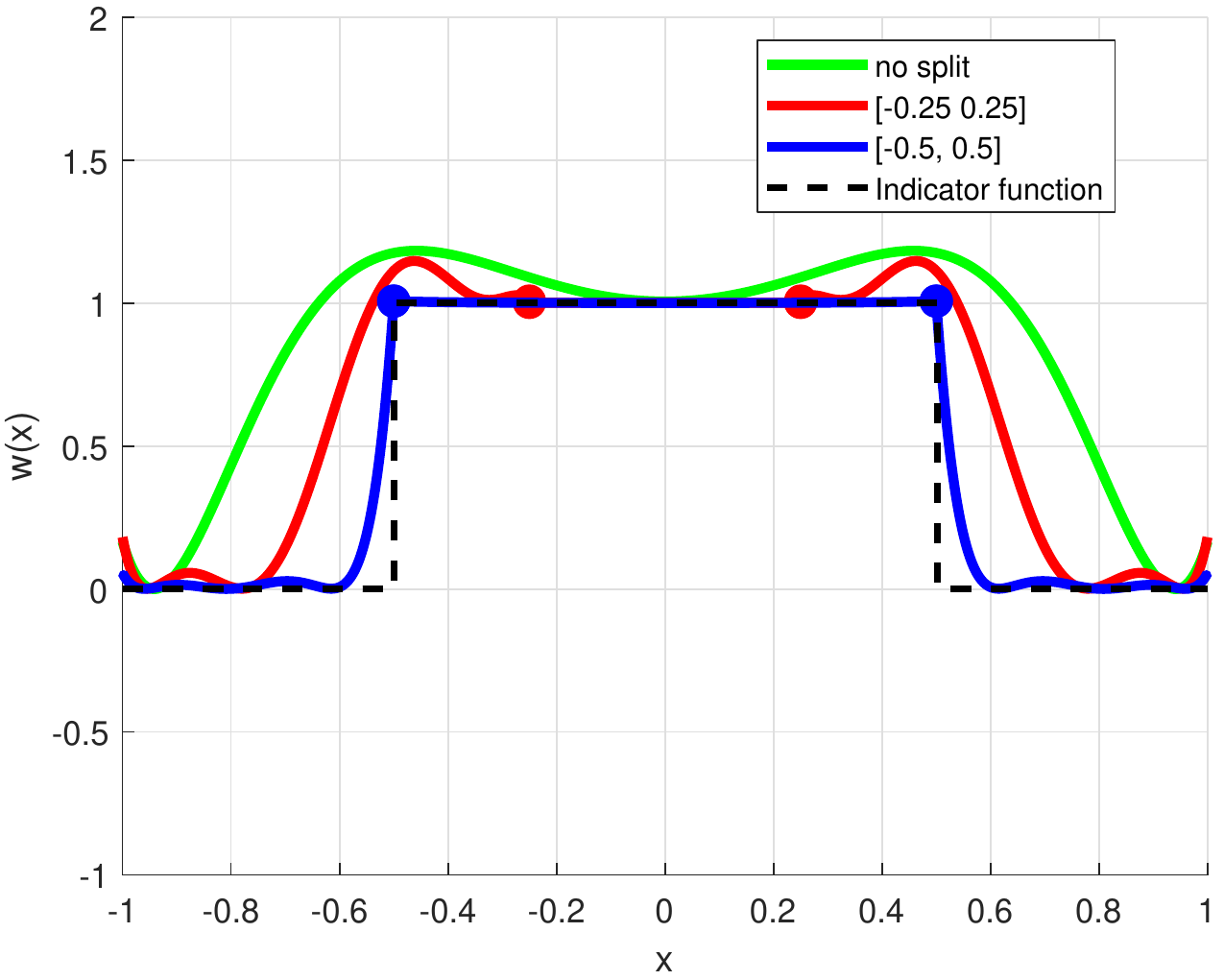}
	\caption{Univariate cubic dynamics,
	approximations of the indicator 
	function with degree 8 polynomials and two splits.
	The ROA is given by \mbox{\(\{x:w(x) \geq 1\}\)}.
	The approximation is more precise
	for splits that are close to the bounds of the ROA (red), 
	and gives the exact ROA when the splits are exactly at the bounds (blue).
	Both split polynomials are more precise than the green, non-split one.}
	\label{fig:cubic_low}
\end{figure}

In Fig. \ref{fig:cubic_low}, we can see a comparison between 
the original method (without splits) and multiple calculations with splits,
going from the inside of the real ROA to the outside.
The ROA estimates here are given by \( \{x:w(x) \geq 1\}\),
which follows from \eqref{eq:def_X0v} and \eqref{eq:dualOriginalSPlit}.
We can observe that the ROA estimates get more precise, the closer the 
splits are to the real ROA.
Let us define the estimate of \(I_{X_0}(x)\) as \(\bar{I}_{X_0}(x)\), which takes \(1\) on	\(\{x:w(x)\geq 1\}\) and \(0\) otherwise.
We can observe from Fig. \ref{fig:cubic_low},
that for the exact split it holds that \mbox{\( \bar{I}_{X_0} = I_{X_0}\)} and we obtain the theoretically optimal estimate.
This example shows that our method can be used in an iterative manner 
with splits along an inner approximation of the ROA (such as the one in \cite{InnerKorda})
as a starting point.

\subsection{Brockett integrator}
\label{sec:brockett}
The Brockett integrator is defined according to 
\cite{Brockett83asymptoticstability} as
\begin{equation}
\begin{alignedat}{2}
\dot{x}_1 &= u_1\\
\dot{x}_2 &= u_2\\
\dot{x}_3 &= u_1x_2 - u_2x_1.
\end{alignedat}
\end{equation}
With \(X = \{x \in \mathbb{R}^3 : ||x||_\infty \leq 1\}\), \(X_T = \{0\}\), \mbox{\(U = \{u \in \mathbb{R}^2 : ||u||_2 \leq 1\}\)}, and \(T = 1\).
As was stated before, this system usually serves as a benchmark 
for nonholonomic control strategies, because it is the simplest system for which 
there exists no continuous control law which would make the origin asymptotically
stable \cite{Brockett83asymptoticstability}.

We shall use the system for calculation of the controlled ROA,
which can be computed analytically \cite{Lasserre2007} as
\begin{equation}
\begin{alignedat}{2}
\mathcal{T}(x) = \frac{\theta\sqrt{x_1^2 + x_2^2 + 2|x_3|}}{\sqrt{\theta + \sin^2 \theta - \sin\theta \cos \theta }},
\end{alignedat}
\end{equation}
where \(\theta = \theta(x)\) is the unique solution in \([0,\pi)\) to 
\begin{equation}
\begin{alignedat}{2}
\frac{\theta - \sin \theta \cos \theta}{\sin^2 \theta} (x_1^2 + x_2^2) = 2|x_3|.
\end{alignedat}
\end{equation}

\begin{figure}[h]
	\centering
	\includegraphics[trim=0 0 0 -1.13cm, width=0.5\textwidth]{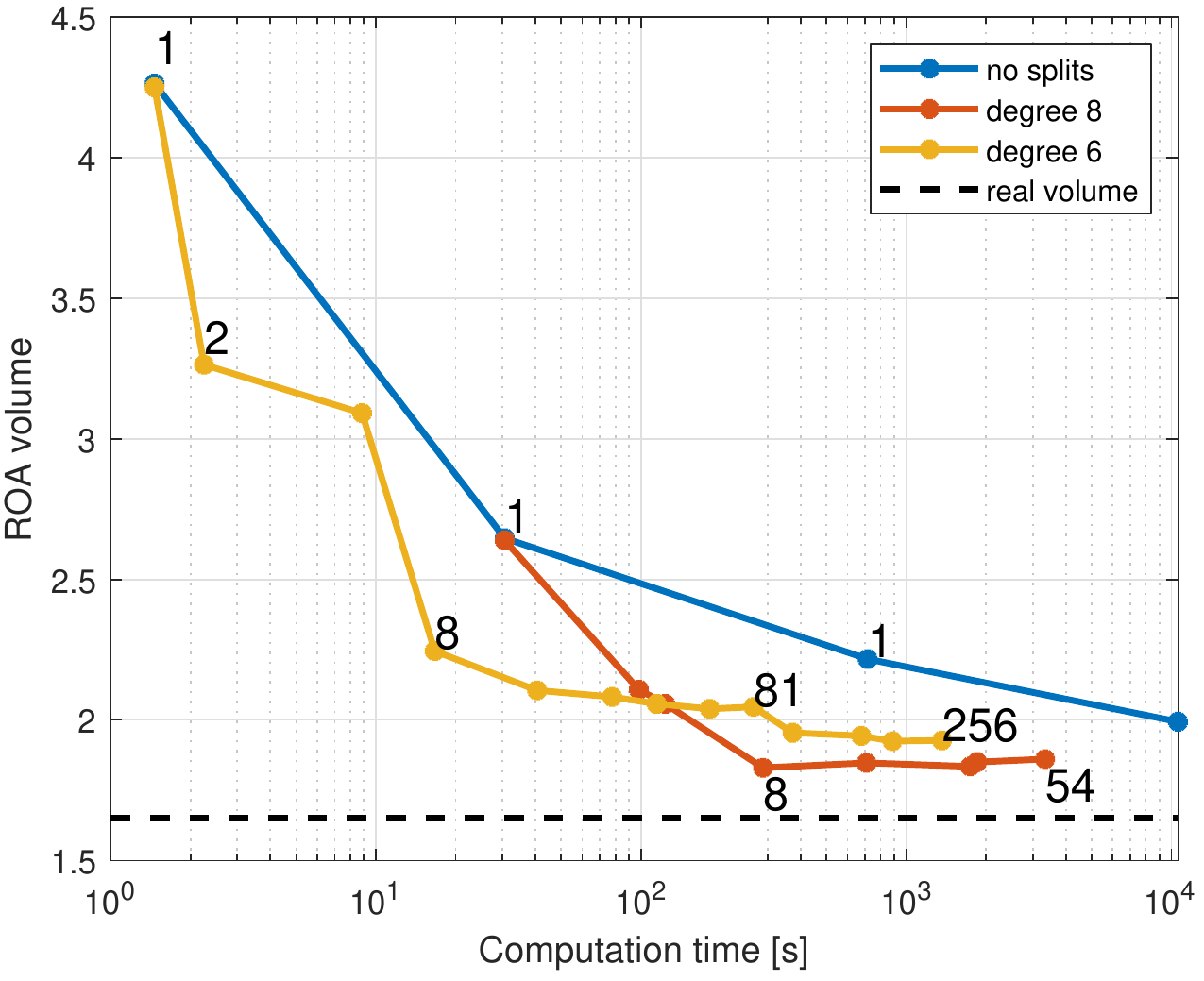}
	\caption{Comparison of various methods on the Brockett integrator.
		The blue line shows the original algorithm without splits and with increasing degree
		of the approximation \(d_\text{o} \in \{6,8,10,12\}\).
			The other lines show the performance of the proposed approach, each having 
			a fixed degree \mbox{\(d_\text{s} \in \{6,8\}\)} with increasing numbers of cells \(X_i\), which are
			denoted as numbers next to the datapoints. The volumes were estimated 
			using Monte-Carlo methods.}
	\label{fig:bro_comptime_vs_volume}
\end{figure}

The Fig. \ref{fig:bro_comptime_vs_volume} shows that given a fixed 
time for the calculation, the proposed approach is always better than the 
original one and that the split-method approaches the real volume much faster,
although neither of the two methods reached the real volume, due to memory constraints.

\begin{figure}[tbhp]
	\centering
	\includegraphics[width=0.5\textwidth]{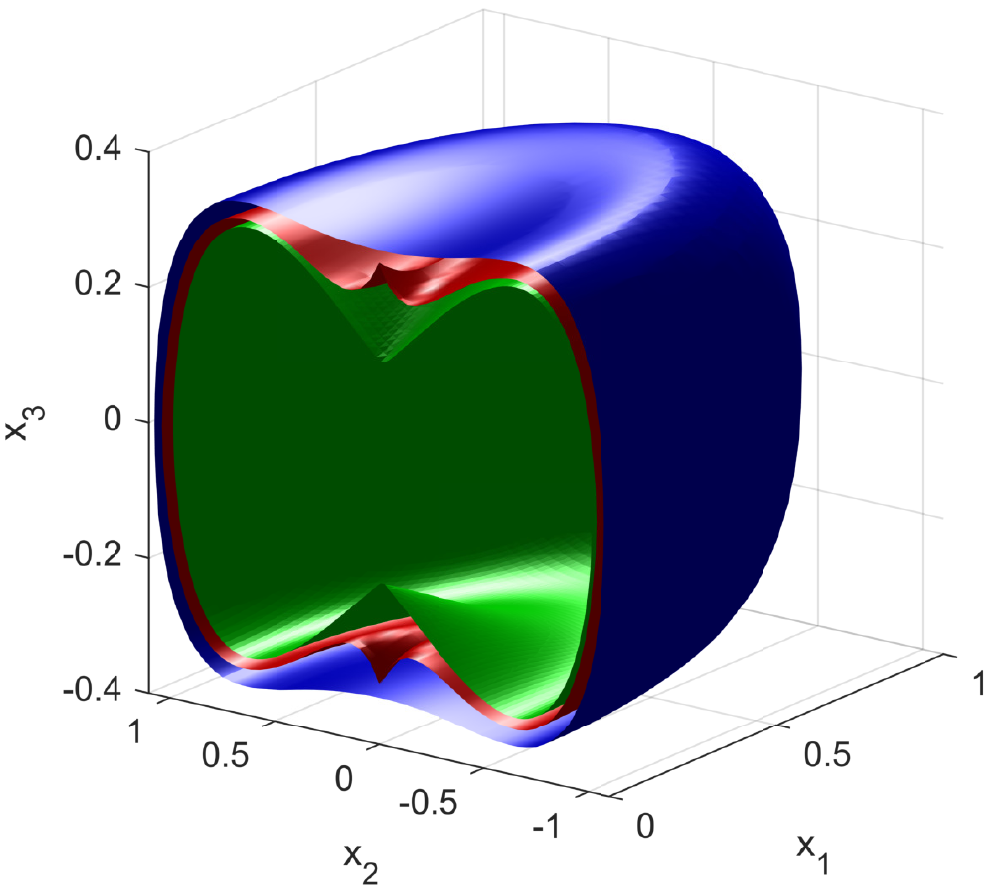}
	\caption{Sliced ROA of the Brockett integrator (green) and two slices of its 
	approximations with similar 
		computation time.
	Blue approximation is by single degree 10 polynomial (754s)
	and the red is by sixteen polynomials of degree~8 (753s).
	The red, lower-degree, approximation is visibly closer to the real ROA.
	}
	\label{fig:bro}
\end{figure}

A visual example of the difference between the two methods can be seen 
in the Fig. \ref{fig:bro} where two ROA's with the same computation time 
are compared, with one being calculated by the original method and the other 
by the proposed method. We can see that there is a notable difference between the 
two approximations, 
and that the better approximation is done by the lower-degree polynomials.

\subsection{Performance and scalability}
\label{sec:performance}
\subsubsection{Problem size}
First, we shall investigate the accuracy of the algorithm with increasing size of the SDP.
The problem size is measured as the number of nonzero elements in the 
\(A\) matrix of the SDP
\begin{equation}
\begin{alignedat}{2}
&\min \quad c^{\top} x \\
&\text{s.t. } Ax = b, \quad x \in \mathcal{K}
\end{alignedat}
\end{equation}
for variable \(x \in \mathbb{R}^n\), convex cone \(\mathcal{K}\) and
data  \(A \in \mathbb{R}^{m \times n}\), \(b \in \mathbb{R}^m\), and \(c \in \mathbb{R}^n\).

\begin{figure}[htbp]
	\centering
	\includegraphics[width=0.5\textwidth]{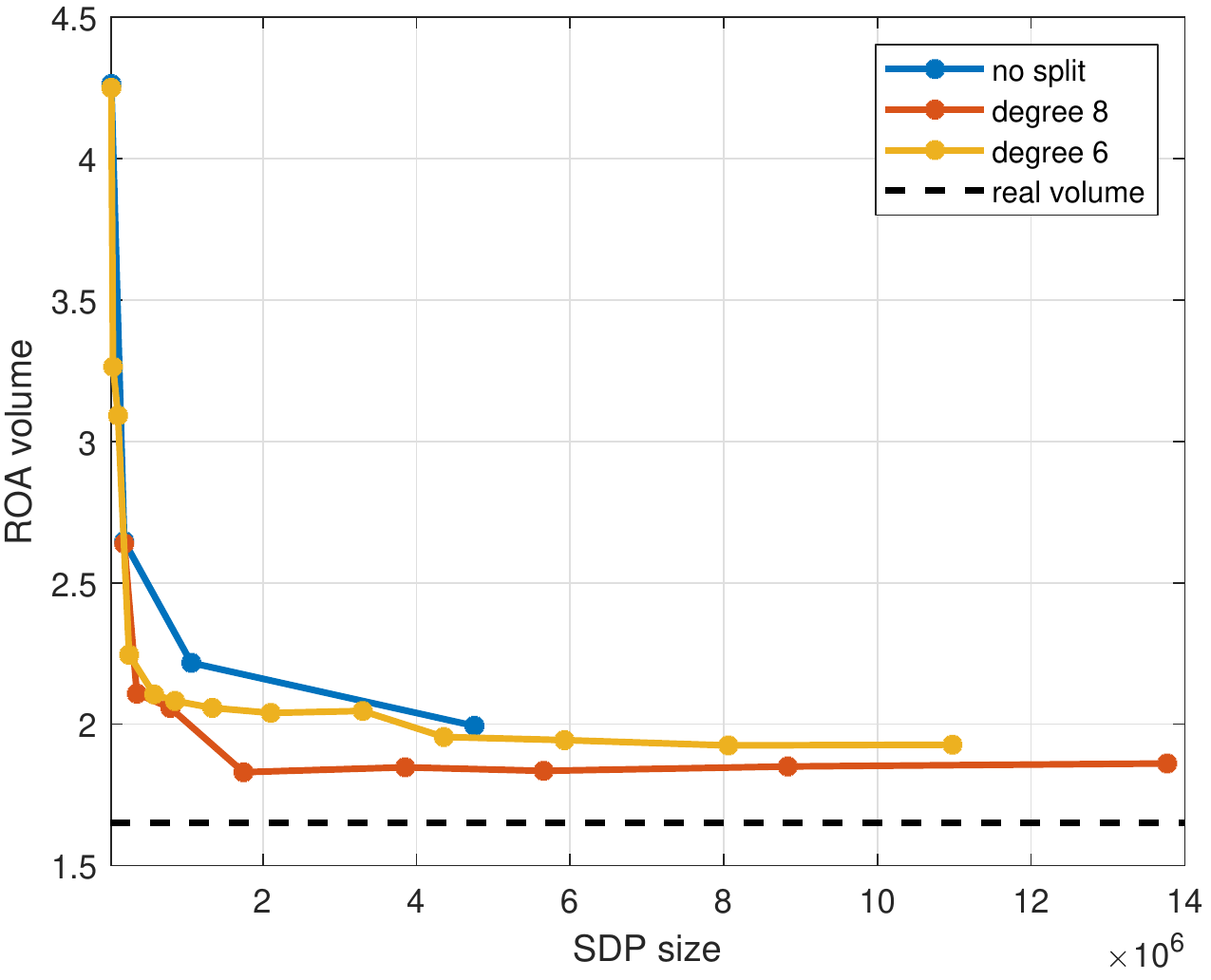}
	\caption{Brockett integrator - Given SDP size, the split problem always give better results. The volumes were estimated 
			using Monte-Carlo methods.}
	\label{fig:nnz_volume}
\end{figure}

We can clearly see in Fig. \ref{fig:nnz_volume} that
the split versions are always more precise than the non-split version
with the same memory footprint.

\subsubsection{Computation time}
\begin{figure}[htbp]
	\centering
	\includegraphics[width=0.5\textwidth]{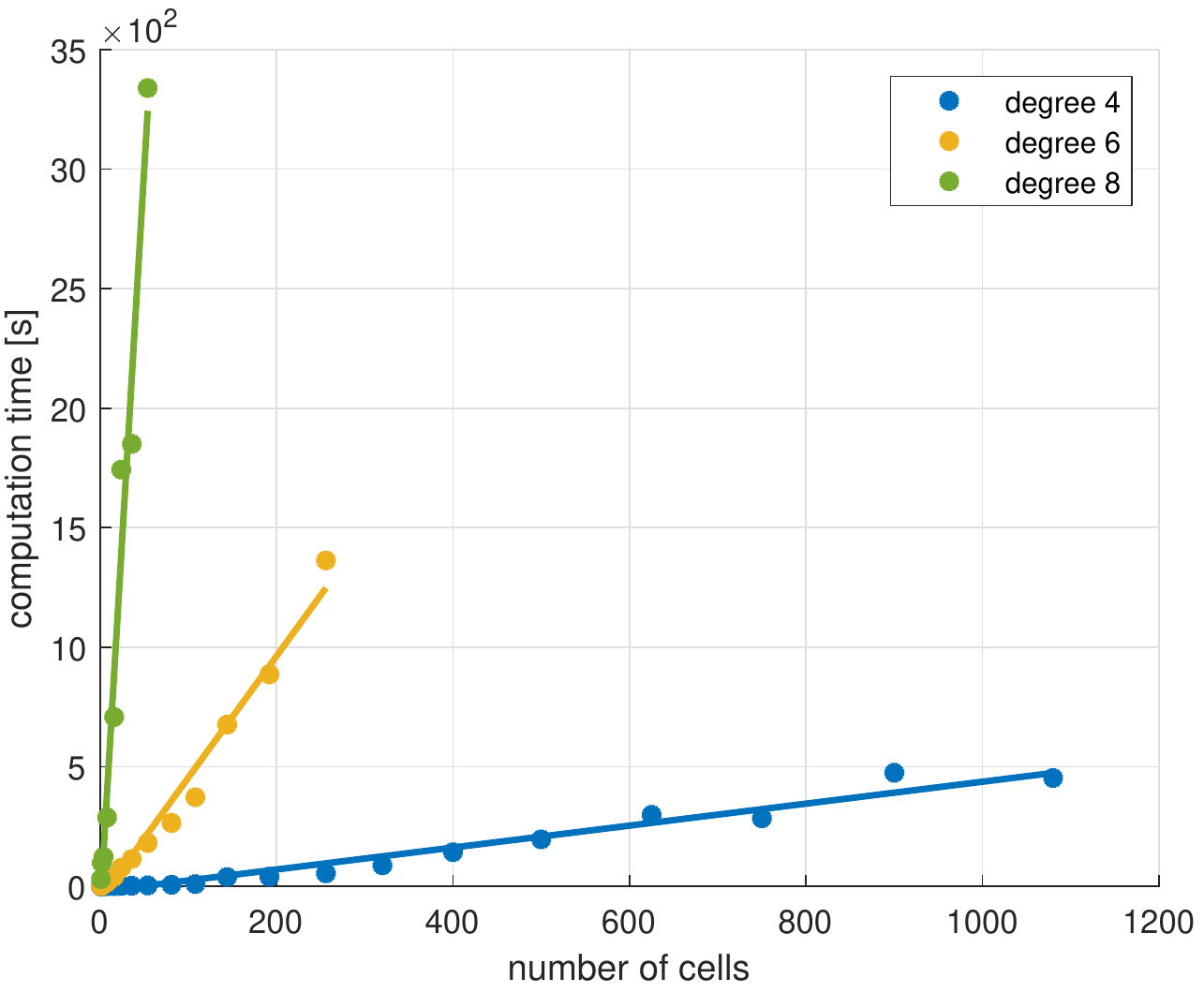}
	\caption{Brockett integrator - The computation time increases linearly 
	with the number of cells \(X_i\).}
	\label{fig:bro_subspace_comptime}
\end{figure}

\begin{figure}[htbp]
	\centering
	\includegraphics[trim = 0 0 0 -0.355cm , width=0.5\textwidth]{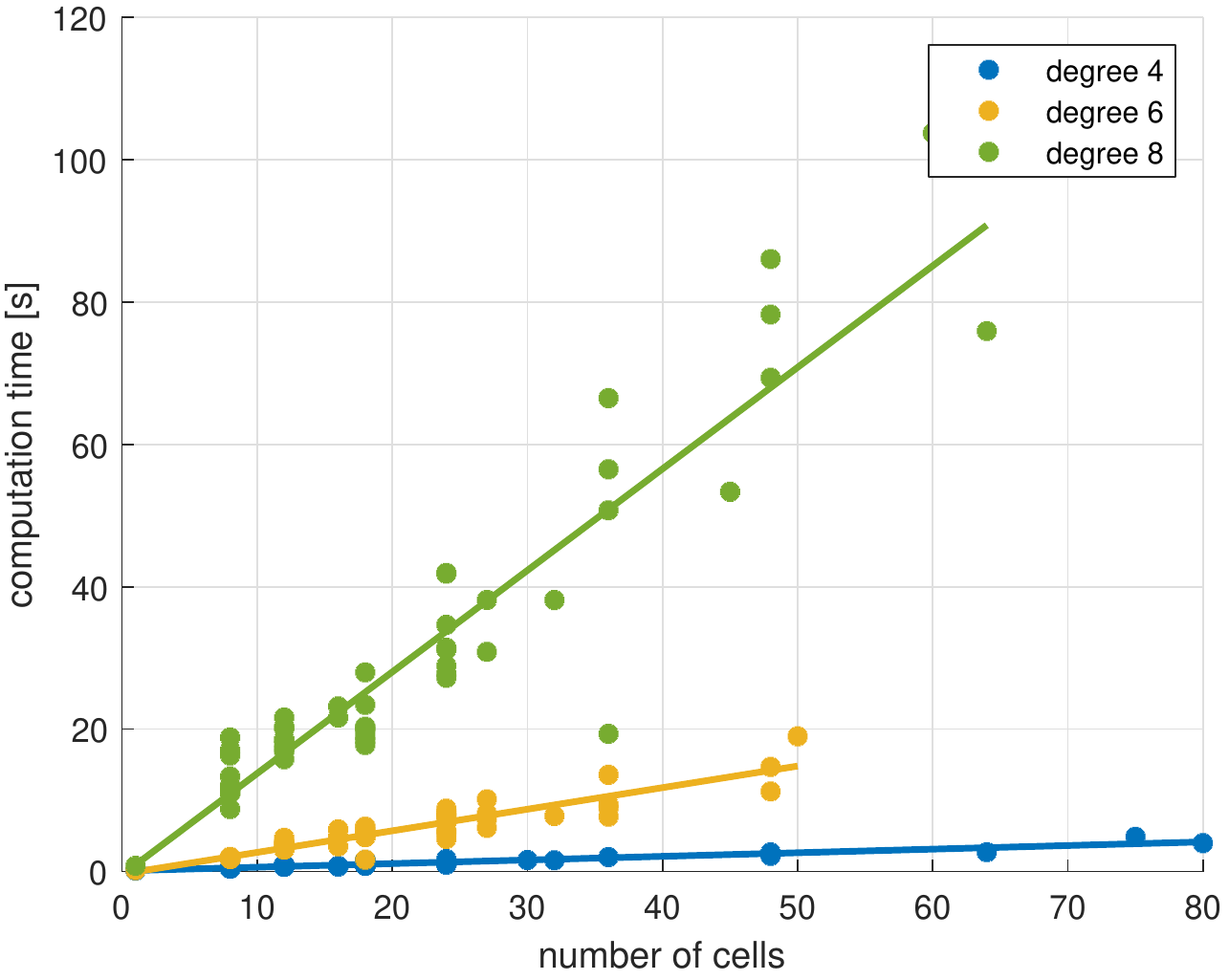}
	\caption{Double integrator - The computation time grows 
	linearly with the number of randomly chosen cells \(X_i\).}
	\label{fig:dual_subspace_comptime}
\end{figure}

The problem size increases linearly with the number of cells (as was explained in \ref{sec:practical_impl}), 
but the computation time 
does not necessarily have to follow the same pattern.
In this case, however, the computation time also showed
linear growth as can be seen in the Figures \ref{fig:bro_subspace_comptime}
and \ref{fig:dual_subspace_comptime} for the Brockett integrator and the Double
integrator respectively. 
\begin{samepage}The Double integrator is defined as
 \[\dot{x}_1 = x_2, \dot{x}_2 = u\] with 
\(X = [-0.7\times 0.7]\times[-1.2 \times 1.2]\), \(X_T = \{0\}\) and \(T = 1\).
See \cite[9.3]{milanroa} for more details.
\end{samepage}
In the Fig. \ref{fig:dual_subspace_comptime}, the variables to be split were 
chosen randomly and the splits were always halving the largest interval of 
the randomly selected variable.
This was done in order to ensure that the linear growth is not
 simply a fortunate result of a particular split order.


\addtolength{\textheight}{-3cm}   

\section{Conclusion}
\label{sec:conclusion}
We showed that the convex SOS-based calculation of ROA can be 
extended by splitting the time and state space to achieve
better accuracy with lower degree polynomials, while keeping 
the outer-approximation guarantees.
The method is faster and provides
more accurate estimates of the ROA. 
It was also empirically demonstrated that the computation time grows linearly with increasing
number of cells.

\section{Future work}
\label{sec:future_work}
The sparse problem structure can be exploited by ADMM-like algorithm,
where the ADMM could alternate between solving small independent SDPs, and 
enforcing the flow constraints between them. This would allow for parallelization and significantly decrease the required computation time.

The splits can be found in an optimization-based fashion by differentiating the SDP \cite{diff_sdp} with respect to the location (or more general parametrization of the splits), thereby increasing the accuracy while keeping the memory requirements constant. This also leads to the option of creating the splits in a less naive
fashion than the one presented here. Using inner-approximation of the ROA 
lends itself as an obvious candidate, 
due to the behaviour shown in \ref{sec:univar}.

The same splitting approach can also be applied to other SOS-based algorithms,
e.g., to the optimal control problem \cite{optimal_control_LMI}, where the increased
precision would lead to more precise estimates of the optimal control sequences
and also allow for control of higher dimensional systems. When applied to the problem of the maximum control invariant set estimation~\cite{korda2014mci}, temporal splitting is no longer possible but could be replaced by introducing multiple discount factors.

Finally, a natural way to decrease the computational complexity further is to combine the approach with the sparsity-based decomposition method of~\cite{schlosser2020sparse}, where the time and state-space splitting would be applied to the subsystems that cannot decomposed using the method of~\cite{schlosser2020sparse}.







\bibliographystyle{IEEEtran}
\bibliography{IEEEabrv,./library_roa_split.bib}

\end{document}